# Inventory Management Under Stochastic Demand: A Simulation-Optimization Approach

Sarit Maitra
Alliance Business School, Alliance University
Sarit.maitra@gmail.com

**Abstract**

This study presents a comprehensive approach to optimizing inventory management under stochastic demand by leveraging Monte Carlo Simulation (MCS) with grid search and Bayesian optimization. By using a business case of historical demand data and through the comparison of periodic review (p, Q) and continuous review (r, Q) inventory policies, it demonstrates that the (r, Q) policy significantly increases expected profit by dynamically managing inventory levels based on daily demand and lead time considerations. The integration of random and conditional sampling techniques highlights critical periods of high demand, providing deeper insights into demand patterns. While conditional sampling reduces execution time, it yields slightly lower profits compared to random sampling. Though Bayesian optimization marginally outperforms grid search in identifying optimal reorder quantities and points, however, given the stochastic nature of the algorithm, this can change with multiple runs. This study accentuates the effectiveness of advanced simulation and optimization techniques in addressing complex inventory challenges, ultimately supporting more informed and profitable inventory management decisions. The simulation model and optimization framework are open-source and written in Python, promoting transparency and enabling other researchers and practitioners to replicate and build upon this work. This contributes to the advancement of knowledge and the development of more effective inventory management solutions.

## 1. Introduction

According to the 34th Annual Council of Supply Chain Management Professionals State of Logistics Report 2023, U.S. business logistics costs stood at a record $2.3 trillion (about $7,100 per person in the US) (about $7,100 per person in the US, which is around 9% of US GDP). The criticality of logistics for an organization's business processes and supply chain network is prevalent in the empirical studies (see [1]; [2]). While the traditional supply chain used to be modest and simple, the global economy has driven it to evolve and expand into an enormous network of chains that extend outside of a nation's geographical bounds. Consequently, inventory management has evolved into a more strategic function since it helps businesses establish and maintain a competitive advantage [3].

Empirical studies ([4]; [5]) highlighted the four major areas of logistical costs that many businesses have a blind spot. These are (a) inventory management costs; (b) warehouse costs; (c) transportation costs; and (d) distribution costs. There is no denying that effective inventory management has a direct impact on organizational performance (see [6]; [7]). However, inventory management is a complex task, and the intricacies are accentuated when confronted with stochastic demand or demand uncertainty, which leads to stockouts, excessive inventory holdings, revenue losses, etc. In the context of tactical supply chain management, traditional operations research approaches continue to confront significant challenges in practice [8].

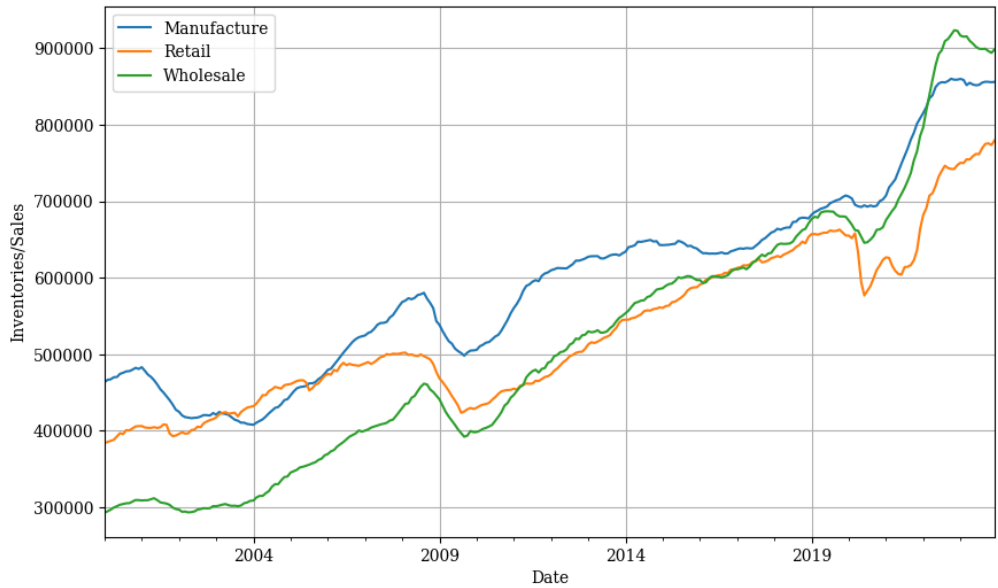

**Figure 1.** Seasonally adjusted inventories/sales (Source: U.S. Census Bureau, 2023)

**Figure 1** displays a plot generated from the data sourced from the US Census Bureau (2023)[1] showing an overall increase in inventories and sales from 2000 to 2023 in three sectors (wholesale, retail, and manufacturing). The plot reflects the importance of inventory for economic growth, and thus optimization is a serious activity to improve service levels and inventory policy [9]. In the supply chain, optimization is essential for exception management which ensures steady operations despite potential disruptions and is crucial for competitiveness.

Despite being computationally efficient, optimization using analytical models is impractical in supply chain settings [10] due to the uncertain or stochastic business dynamic. Here simulation models offer the flexibility to incorporate any stochastic components and typically enable the modeling of all the intricacies and dynamics of real-life scenarios without making too simplistic assumptions. In recent times, researchers have turned their focus to stochastic demand modeling (see [11]; [12]; [13]; [14]; [15]). In that context, Pinder [16] cautioned that decision-making under uncertainty could be highly inaccurate and have negative effects on the business if simulation is not applied. Stochastic demand inherently involves uncertainty and variability (see [17]; [18]), which makes it a probabilistic model.

The primary focus of an inventory model is to determine how much to order and when to order, which historically has been focused on deterministic demand models (see [19]; [20]; [21];[22]; [23]). The applicability of such a model in real-world scenarios has been increasingly limited due to a long list of constraints and the rising uncertainty in demand patterns. In line with recent advancements, Ekren et al. [24] have introduced a demand forecasting model that not only considers demand uncertainty but also employs a robust optimization approach. Other researchers, too, made clear arguments on optimization techniques, which play an indispensable role in implementing an effective inventory system (see [25]; [26]; [27]). To explore nearly optimal control parameters for a stochastic multi-product inventory control system, Jackson [28] suggested a simulation-optimization framework, allowing for full inventory dynamics with risk and reliability analysis. More empirical studies have highlighted the importance of simulation modeling, which helps decision-makers analyze different scenarios and choose the best course of action (see [29]; [30]; [9]; [31]; [32]; [33]; [34]; [35]).

---

[1] https://www.census.gov/mtis/current/index.html

To put this into perspective, simulation modeling simulates the functioning of an existing system and tests various scenarios to provide evidence for decision-making. The goal of inventory optimization is to maximize the performance of the objective function ($f(x)$) in a stochastic environment by identifying an optimal set of control parameters for inventory strategy [36]. The significance of adopting data-driven, analytical, and systematic approaches is imperative to deal with the complex challenges in modern inventory management systems [1]. While several empirical works (e.g., [37]; [38]; [39], among others) delved into the Economic Order Quantity (EOQ) approach, they lacked rigorous analytical evaluations to substantiate their findings. In their work, Homem-de-Mello & Bayraksan [40] highlighted that the stochastic problem bridges two distinct research communities: optimization and simulation. The optimization community involves structural exploitation, theoretical foundations, and algorithm development, whereas the simulation community involves a black-box approach with probabilistic modeling. This study combines the strengths of both communities to develop a hybrid method that is powerful and versatile. The proposed approach allows for an inclusive exploration of the solution space to enhance the effectiveness of decision-making. To establish the strategy, Monte Carlo simulation (MCS) is used to generate demand and lead time distributions for various products based on historical demands and estimate potential profit maximization under different inventory policies.

The study develops a model based on MCS and employs (p, Q) and (r, Q) inventory policies. The model is assessed using a case study of different product sales and seeks to implement the right inventory strategy to manage its inventory. The goal is to maximize the predicted profit. The study further uses a grid search between a range of values for each product, a combination of order quantity and re-order points to optimize the profit function. The finding shows that (r, Q) significantly improves the profitability compared to (p, Q) policy, and optimization further improves the profit by 18.64%. The system further employs an optimization approach using both Grid-search and Bayesian optimization along with different inventory management metrics such as average profit ($\mu_{profit}$), reorder quantity (r), order quantity (Q), average safety stock ($\mu_{SS}$), and standard deviation of profit ($\sigma_{profit}$). Bayesian optimization appears to offer marginal improvement, but potentially riskier inventory management strategy compared to grid-search. This trade-off between profit and risk should be considered when deciding which optimization method to use.

To summarize, this study combines MCS and optimization techniques that utilize grid search and Bayesian optimization. The study provides empirical evidence through case studies to validate the model. These contributions highlight the study's originality and its potential impact on improving inventory management practices in supply chain management.

## 1.1 Research gap and study contribution

Empirical work with a case study on stochastic demand in the context of inventory management is limited [41]. Moreover, though simulation is a popular approach to address practical business issues, there are still certain unresolved research issues with simulation-based optimization techniques in the current study context. Simulations are typically computationally expensive since several iterations are needed to overcome noise in the returned result. Furthermore, simulation provides a "what-if" reaction to system inputs in the absence of available gradient information. While specialized optimization techniques such as Simulated Annealing (SA) and Genetic Algorithms (GA) can be used, empirical studies have highlighted that they are all metaheuristic approaches that cannot ensure the quality of the solution (see [42]). Though surrogate models can be adapted for efficient and effective exploration of complex design spaces [43] however, inventory optimization presents unique challenges that complicate the application of surrogate models. Since a supply chain network is typically made up of several facility nodes, and if numerous control parameters need to be calculated, the problem's dimension may be quite significant. Moreover, the inventory management problem's objective and constraints are both opaque, making the task of creating the surrogate function for each equation extremely difficult. Furthermore, a lot of the surrogate-based techniques currently in use are not very broad and rely on a lot of premises, including divergent networks, normally distributed demand, convexity assumptions, etc. Therefore, the current study delves into the development of an effective easy-to-use simulation optimization algorithm to deal with more generic inventory optimization issues under uncertainty. By prioritizing profit maximization, the study aligns more closely with the primary business objectives of many businesses, offering a practical tool for decision-makers who aim to enhance their financial performance through inventory and pricing synergy.

The two main contributions of this paper are:

- This work presents a simulation-optimization model designed to determine near-optimal inventory policies for products with demand uncertainties. This model is adaptable to various probability distributions, thus accommodating uncertainties without imposing strict distribution restrictions.
- Bayesian optimization appears to offer marginal improvement, but potentially riskier inventory management strategy compared to grid-search. This trade-off between profit and risk should be considered when deciding which optimization method to use.

Moreover, building upon the empirical evidence on the efficacy of inventory systems ([44]), this study contributes to the existing body of knowledge by offering a comprehensive examination of popular inventory management models with a case study.

## 2. Methodological framework

**Figure 2** displays a simplified framework for the study's theoretical framework. This framework is further improvised to meet this study's goal and the complexity of the work. In the final stage, the optimal inventory policy was determined which maximizes the profit considering the demand scenarios.

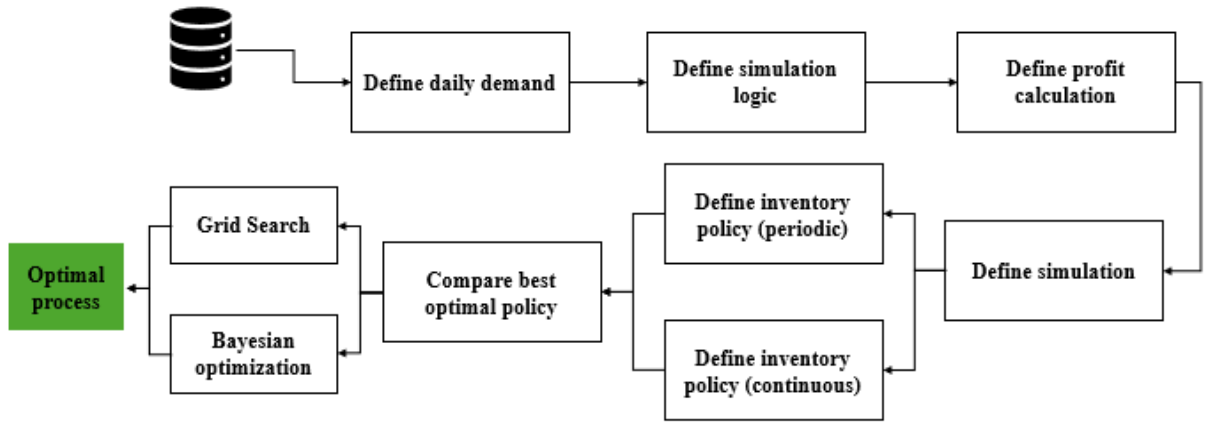

**Figure 2.** Simplified framework.

The following notations are used to formulate the algorithm, which captures the daily inventory dynamics, ordering policy, and profit calculation over 365 days. The profit calculation here aligns with the operational goals, by adjusting holding and ordering costs based on data.

| Notation | Description |
|---|---|
| $D_t$ | Demand on day t |
| $I_t$ | Inventory level at the end of the day t |
| $S_t$ | Units sold on day t |
| $R$ | Review period (number of days) |
| $Q$ | Order quantity placed at each review period |
| $LT$ | Lead time (days) for replenishment |
| $P_s$ | Selling price per unit |
| $C_P$ | Purchase cost per unit |
| $C_H$ | Holding cost rate per unit per day |
| $C_O$ | Ordering cost per order |
| $n$ | Number of review periods in a year $\left(n = \frac{365}{R}\right)$ |
| $SS$ | Safety stock |

| $\delta$ | Indicator variable that takes the value 1 if an order is placed, and 0 otherwise |
|---|---|

For each day $t$ ($t$ ranges from 1 to 365), the inventory and sales are updated as in Eq. (1).

$$S_t = min\ (I_{t-1}, D_t) \qquad (1)$$

Here, $min\ (I_{t-1}, D_t) \rightarrow$ the number of units sold on day $t$ ($S_t$), and this is based on available inventory from the previous day ($I_{t-1}$) and the demand on the current day ($D_t$). The logic behind this is that if the demand $D_t \leq I_{t-1}$, then $S_t = D_t$ which means all the demand can be satisfied. If the demand $D_t \geq I_{t-1}$, then, $S_t = I_{t-1}$ which means only as much inventory as is available can be sold. Both the demand on the current day and the existing inventory from the previous day limit the number of units sold on any given day (t). This ensures that sales never exceed the available inventory, and that demand is fulfilled as much as possible given the constraints. After fulfilling the demand for the day, the remaining inventory is updated by subtracting the units sold from the inventory level of the previous day Eq. (2).

$$I_t = I_{t-1} - S_t \qquad (2)$$

Orders are placed at regular intervals which is the review period (p). If the current day $t$ is a multiple of p (i.e., $t\ mod\ p = 0$), an order is placed to replenish the inventory. Eq. (3) presents the immediate inventory update where for simplicity it is assumed that the inventory is updated immediately with the reorder quantity ($Q$).

$$I_t = I_t + Q \qquad (3)$$

However, in a realistic scenario, orders placed do not arrive instantly but after a lead time ($LT$). Hence, the inventory is updated after this lead time as shown in Eq. (4).

$$I_{t_0+LT} = I_{t_0+LT} + Q\ (if\ t_0 + LT \leq 365) \qquad (4)$$

Where, $t_0$= time or day when the order was placed, $t_0 + LT$ = day when the order arrived.

With the above conditions, the proposed system will simulate inventory levels, sales, and reorder quantity over the 365 days, to calculate the profit based on sales revenue and costs (holding, ordering, and purchase costs). Eq. (5), (6), (7), and (8) formulate the revenue, holding costs, ordering costs and purchase costs.

$$revenue = P_s \sum_{t=1}^{365} S_t = \Sigma(units_{sold}) * (selling\ price) \qquad (5)$$

$$holding\ cost = C_H \sum_{t=1}^{365} I_t = \frac{\Sigma(inventory\ level) * (holding\ cost) * (size)}{days} = \frac{\sum_{t=1}^{365} I_t * C_H * size}{365} \qquad (6)$$

$$\text{ordering cost} = C_O = \text{number of orders} \times \text{ordering cost per order} \quad (7)$$

$$\text{purchase cost or cost of goods (COGS)} = C_P \sum_{t \text{ where } t \bmod p = 0} Q \quad (8)$$

$$\text{total costs} = \text{holding cost} + \text{cost of goods sold} + \text{ordering cost}$$

Since costs are related to a specific product and profit is dependent on demand linearly, the maximum for that function can be determined. Eq. **Error! Reference source not found.**) formulates the simplified profit for the entire year.

$$\text{profit} = \text{revenue} - \text{total costs}$$

$$= P_s \sum_{t=1}^{365} S_t - \left\{ \left( C_H \sum_{t=1}^{365} I_t \right) + (C_O * n) + \left( C_P * \sum_{t \text{ where } t \bmod p = 0} Q \right) \right\} \quad (9)$$

Given the review period and constant demand, the average inventory level $\left( \frac{\sum_{t=1}^{365} I_t}{365} \right)$ can be approximated by $\frac{Q}{2} * 365$ as displayed in Eq. (10).

$$\frac{\sum_{t=1}^{365} I_t}{365} \approx \frac{Q}{2} * 365 \quad (10)$$

Hence the profit can be rewritten as Eq. (11)

$$\text{profit} = P_s \sum_{t=1}^{365} S_t - \left\{ \left( C_H \frac{Q}{2} * 365 + C_O \frac{365}{R} + C_P \frac{365}{R} \right) \right\} \quad (11)$$

To find the optimal order quantity, Eq. (12) differentiates the profit function concerning Q.

$$\frac{d(Profit)}{dQ} = -C_H \frac{365}{2} - C_P \frac{365}{R} \quad (12)$$

Since the first derivative gives us the rate of change of profit concerning order quantity, setting it to zero gives the optimal point as shown in Eq. (13)

$$-C_H \frac{365}{2} - C_P \frac{365}{R} = 0 \Leftrightarrow Q = Q^* \quad (13)$$

However, this approach shows the constant rate of change rather than deriving a critical point since the original profit function is linearly dependent on the order quantity in a simplified model. Thus, optimizing requires balancing between ordering and holding costs. The objective function for maximizing profit, considering all costs and revenues, is displayed in Eq. (14).

$$objective\ function\ (f(x)):\ P_s \sum_{t=1}^{365} S_t - \left\{ \left( C_H \sum_{t=1}^{365} I_t \right) + (C_O * n) + \left( C_P * \sum_{t\ where\ t\ modulo\ p=0} Q \right) \right\} \quad (14)$$

Despite the linear dependency of profit on demand, the resulting optimization problem is not well-represented by a simple mathematical model due to significant uncertainties in demands.

## 2.1. Data exploration

To verify the applicability of the model, an empirical evaluation is provided through a case study. The business case selected here examines four highly customized products and thus demands are unique to every customer. These are highly unpredictable demands with seasonality and latent trends. One year (365 days) of sales information was collected for each product. **Table 1** displays a subset of demand data showing the historical demands of the products for 30 days.

Table 1. A glimpse of varied demand was observed for 30 days.

| Days | Pr1 | Pr2 | Pr3 | Pr4 | Days | Pr1 | Pr2 | Pr3 | Pr4 |
|---|---|---|---|---|---|---|---|---|---|
| 1 | 90 | 610 | 204 | 153 | 16 | 50 | 684 | 166 | 0 |
| 2 | 94 | 685 | 244 | 147 | 17 | 0 | 699 | 193 | 0 |
| 3 | 0 | 649 | 197 | 0 | 18 | 94 | 633 | 0 | 0 |
| 4 | 110 | 667 | 0 | 0 | 19 | 126 | 658 | 260 | 144 |
| 5 | 106 | 663 | 228 | 0 | 20 | 16 | 636 | 205 | 0 |
| 6 | 102 | 663 | 201 | 156 | 21 | 86 | 655 | 245 | 0 |
| 7 | 94 | 681 | 215 | 0 | 22 | 102 | 590 | 201 | 156 |
| 8 | 154 | 671 | 170 | 0 | 23 | 0 | 660 | 223 | 0 |
| 9 | 108 | 653 | 219 | 0 | 24 | 84 | 672 | 0 | 0 |
| 10 | 74 | 631 | 189 | 0 | 25 | 124 | 601 | 231 | 147 |
| 11 | 80 | 663 | 238 | 0 | 26 | 196 | 658 | 0 | 0 |
| 12 | 152 | 707 | 231 | 0 | 27 | 102 | 659 | 0 | 0 |
| 13 | 102 | 650 | 0 | 0 | 28 | 0 | 636 | 262 | 0 |
| 14 | 90 | 692 | 165 | 0 | 29 | 66 | 647 | 181 | 150 |
| 15 | 64 | 649 | 208 | 156 | 30 | 106 | 602 | 231 | 147 |

**Table 2** displays the statistical summary of demand data. Some general inferences are made from this statistical summary as stated below.

Table 2. Statistical summary.

|  | Pr1 | Pr2 | Pr3 | Pr4 |
|---|---|---|---|---|
| Purchase Cost ($C_p$) | € 12 | € 7 | € 6 | € 37 |
| Lead Time (LT) | 9 days | 6 days | 15 days | 12 days |
| Size | 0.57 | 0.05 | 0.53 | 1.05 |
| Selling Price ($P_s$) | € 16.10 | € 8.60 | € 10.20 | € 68 |
| Starting Stock ($S_t$) | 2750 | 22500 | 5200 | 1400 |
| Mean demand over lead time ($mu$) | 103.50 | 648.55 | 201.68 | 150.06 |
| Std. Dev. of demand ($\sigma_D$) | 37.32 | 26.45 | 31.08 | 3.21 |
| Probability (Probability of demand on any day) | 0.76 | 1.00 | 0.70 | 0.23 |
| Order Cost ($C_o$) | € 1000 | € 1200 | € 1000 | € 1200 |
| Holding Cost ($C_h$) → 20% of unit cost | 20% | 20% | 20% | 20% |
| Expected demand (lead time) | 705 | 3891 | 2266 | 785 |
| Std. Dev. ($\sigma_{lead\ time}$) | 165.01 | 64.78 | 383.33 | 299.92 |
| Annual demand | 28,670 | 237,370 | 51,831 | 13,056 |

Given the high purchase cost, comparatively large size, and low probability of demand (0.23), it can be concluded that the demand for Pr4 is infrequent. The $\sigma_D$ is low (3.21) which supports the above observation. There are multiple instances where the demands for Pr1, Pr3, and Pr4 are zero on several days. These can be considered rare events in the context of daily operations when these products are not sold at all. This is particularly noticeable for Pr3 and Pr4 (**Table *1*. A glimpse of varied demand was observed for 30 days.Table 1**).

**Figure 3** displays the distribution plots, which provide insights into the central tendency and dispersion of historical demand data for each product. Pr1 (Top Left): the demand distribution is skewed right (positive skew), indicating that while most days have lower demand, there are occasional high-demand days. Pr2 (Top Right): the demand distribution appears to be normally distributed, with a symmetric bell shape around the mean. This indicates a consistent demand pattern with less variability. Pr3 (Bottom Left): the demand distribution is bimodal, indicating two peaks in demand. This suggests that there are two different demand patterns or seasons affecting this product. Pr4 (Bottom Right): the demand distribution is highly skewed to the right with a significant portion of the demand being zero or close to zero. This suggests that the product is not in demand on many days but has occasional spikes.

Pr1 and Pr3 show higher variability and unpredictability in demand compared to Products 2 and 4. Pr1 and Pr4 have skewed demand distributions, suggesting that most of the demand is concentrated on the lower end with occasional high-demand days. Pr2 has a more consistent and predictable demand, which is closer to a normal distribution.

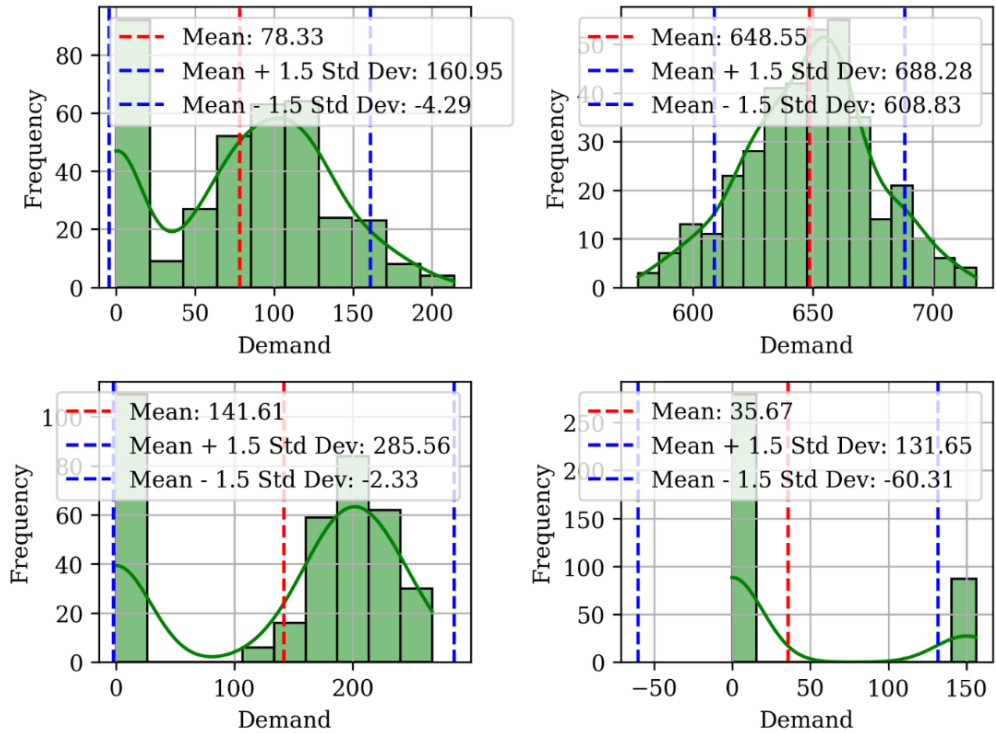

**Figure 3.** Distribution of product demand.

Other than Pr2, no other product was purchased every day. The likelihood that the products were purchased on the day is indicated by the estimated proportion of orders. For example, **Table 2**'s predicted proportion for Pr1 is 0.76, meaning that there was a 76% likelihood that this product would be demanded on any given day of the year.

The inventory levels for each day of the year were added together to display the total amount of stock held over the year. The holdings costs were then calculated by multiplying the quantity of stock owned, the size of the product's unit, and the daily cost of maintaining a unit. **Table 3** displays the pseudocode for the profit calculation.

**Table 3.** Pseudocode for profit calculation.

1. Profit calculation:
2.    Unit Cost, Selling Price, Holding Cost, Order Cost, Size
3.    days = 365
4.    $revenue = P_s \sum_{t=1}^{365} S_t$
5.    $C_o = number\ of\ orders \times ordering\ cost\ per\ order$
6.    $C_h = \frac{\sum(inventory\ level) * holding\ cost * size}{days}$
7.    $profit = \{P_s \sum_{t=1}^{365} S_t - (total\ costs)\}$
8.    RETURN profit.

## 2.2. Bayesian inference and Monte Carlo method

A combination of Bayesian inference principles (probabilistic modeling and inference) and the Monte Carlo method (for numerical estimation) is presented below. The goal is to simulate the daily demand of every product. The major source of uncertainty is the unpredictable demand pattern displayed in Eq. (15) and Eq. (16)

$$\mu \ (Average \ daily \ demand) = \frac{number \ of \ orders \ last \ year}{number \ of \ working \ days} \quad (15)$$

$$variability \ (\sigma) = standard \ deviation \ in \ order \quad (16)$$

MCS generates random samples for demand based on the probability distribution function (PDF) using the mean ($\mu_D$) and standard deviation ($\sigma_D$) of historical demand. Here we introduce conditional sampling and its mathematical formulation. Assuming we have a random variable $X$ representing daily demand. The distribution of $X$ is normally distributed with parameters $\mu$ and $\sigma$. During the simulation, $X$ was conditionally sampled based on whether an order was placed on the previous day.

The probability density function (PDF) of X is formulated as in Eq. (17).

$$f(x) = \frac{1}{\sqrt{2\pi\sigma^2}} exp\left(-\frac{(x-\mu)^2}{2\sigma^2}\right) \quad (17)$$

Here, a condition C is introduced, which indicates whether an order was placed:

- If C = 1: Use demand from the previous day ($X_{prev}$)
- If C = 0: sample X from the usual distribution

For simplicity, we assume

- $f(x|C = 1)$ is a Dirac delta ($\delta$) function
- $f(x|C = 0)$ is a normal distribution of $f(x)$

Eq. (18) presents the integral representation considering these conditions.

$$E(x|C) = \int_{-\infty}^{\infty} x * f(x|C) dx \quad (18)$$

For C = 1:

$$E(X|C = 1) = \int_{-\infty}^{\infty} x * \delta(x - X_{prev}) dx$$

$$E(X|C = 1) = X_{prev}$$

For C = 0:

$$E(X|C = 0) = \int_{-\infty}^{\infty} x * f(x)dx$$

$$E(X|C = 0) = \int_{-\infty}^{\infty} x * \frac{1}{\sqrt{2\pi\sigma^2}} \exp\left(-\frac{(x-\mu)^2}{2\sigma^2}\right) dx$$

In some cases, weighted sampling can be used to correct the bias introduced by conditional sampling. For instance, if we want to correct for the effect of conditioning on C = 1, we might introduce a weight $w(x)$ in the sampling process:

$$E(X|C) \int_{-\infty}^{\infty} x * w(x) * f(x|C)dx$$

where $w(x)$ is a weighting function that adjusts the contribution of each sample $x$ based on its likelihood under the conditional distribution $f(x | C)$. We experimented with conditional sampling to draw a comparison with random sampling. By using conditional sampling, inventory managers can simulate demand scenarios to capture current conditions better. Such as, conditionally sampling demand based on recent historical data or inventory levels allows for modeling how recent trends affect future demands. By adjusting sampling based on inventory levels and expected deliveries, the simulation can predict and manage risks of stockouts or excess inventory. Whether this approach leads to improved profitability, is to be experimented on. Theoretically, integrating cost factors (e.g., holding costs, ordering costs) into conditional sampling helps in optimizing inventory policies that minimize costs while meeting service-level objectives.

MCS calculates the appropriate inventory levels by simulating the inventory and demand dynamics and iterating over each day in one year (365 days). It generates random daily demand based on the demand distribution. When a product is bought, the demand follows a lognormal distribution. The goal is that the model would represent the unpredictability and skewness in demand patterns by utilizing a log-normal distribution. This supposition streamlines the modeling procedure and corresponds with empirical data in numerous situations. A uniform distribution between 0 and 1 was used to select a random number to imitate everyday purchase activity.

Safety stock ensures sufficient inventory is available to meet demand during periods of uncertainty or longer-than-expected lead times. During the review period, the algorithm checks if the inventory level has fallen below the safety stock threshold and accordingly calculates the reorder quantity needed to bring the inventory up to the level of safety stock plus expected demand during lead time.

## 2.3. Demand modeling

Demand per period is expressed in Eq. (19)

$$demand\ per\ period \sim \aleph(\mu_d, \sigma_d^2) \tag{19}$$

The standard deviation of demand is expressed in Eq. (20)

$$\sigma_D = \sqrt{\frac{\sum(d_t - \mu_D)^2}{n-1}} \tag{20}$$

where n is = number of (historical) demand observations and $d_t$ is the demand at time t.

Drawing the reference from the previous section, log-normal distribution was chosen to model the demand over a wide range. The parameters for the log-normal distribution are derived from the $\mu_D$ and $\sigma_D$ during the lead time. The demand over time follows a continuous probability distribution. In this context, the integral helps us understand the expected demand over a period, considering the PDF of the demand. The expected demand $E(D)$ over a period displayed in Eq. (21) which is an integral of the demand $D$ multiplied by its probability density function $f(D)$:

$$E(D) = \int_0^\infty D * f(d) * d(D) \tag{21}$$

For a normal distribution, the PDF is given by Eq. (22):

$$f(D) = \frac{1}{\sigma_d\sqrt{2\pi}} exp\left(-\frac{(D-\mu_d)^2}{2\sigma_d^2}\right) \tag{22}$$

Eq. (23) displays the demand.

$$E(D_{LT}) = LT * \mu_d \tag{23}$$

Eq. (24) shows the variance of demand over the lead time is:

$$\sigma_{D_{LT}}^2 = LT * \sigma_d^2 \tag{24}$$

Safety stock is calculated using the safety factor z, the standard deviation of demand during the lead time, and the lead time itself as in Eq. (25).

$$SS = z * \sigma_d = z * \sigma_D * \sqrt{LT} \tag{25}$$

Eq. (26) shows the order-up-to level ($OUP$) which is the level up to which inventory needs to be replenished to cover the expected demand during the review period $R$ and the safety stock.

$$OUP = E(D_R) + SS = R * \mu_d + z * \sigma_D * \sqrt{LT} \tag{26}$$

The order quantity (Q) at each review period is in Eq. (27)

$$Q = \max(0, OUP - \text{current inventory}) \tag{27}$$

Eq. (28) integrates the demand PDF over a continuous period.

$$E(D_R) = \int_0^R D * f(D) * d(D) \tag{28}$$

Given that $D$ follows a normal distribution, this integral provides the average demand over the review period ($r$).

## 2.4. Simulation and analysis

Two separate sets of simulation logic were set up – (1) for (p, Q) policy and (2) for (r, Q) policy. The logic was set up using MCS to optimize the reorder quantities for the set of products. It calculates the expected profit for different reorder quantities (Q) and identifies the optimal Q that maximizes profit. The algorithm calculates the expected profit for each reorder quantity (Q) by simulating inventory operations and evaluating the profit over 1000 iterations. The simulation model and optimization framework are both open-source and written in Python v3.10.1.

For (p, Q) policy, the algorithm places an order to replenish stock if the total number of days equals the review period. This implies that inventory levels are checked, and orders are placed at regular, fixed intervals (the review period). The reorder quantity is determined and added to the inventory after the review period has passed and the lead time has elapsed. The algorithm checks if the current inventory is sufficient to fulfill the day's demand and is independent of whether the review is periodic or continuous. The inventory level is updated based on daily demand, and the revenue is calculated based on units sold, which applies to both periodic and continuous reviews. In (r, Q) policy, the algorithm would continuously monitor inventory levels and place an order whenever the inventory level falls below a predetermined reorder point ($r$). Orders are not tied to fixed intervals but are triggered by inventory levels.

The concept is predicated on the straightforward estimation that satisfying the day's demand requires only considering the present inventory level. This means that the algorithm makes a simple assumption: to fulfill the demand for any given day, it only needs to check if the current inventory is sufficient. The amount of demand reduces the inventory, whereas the number of units sold that day increases revenue by the same amount. The algorithm counts the number of days in the current period and places an order to replenish the stock (reorder quantity) if the total number of days is equal to the review period. The algorithm receives this value as an input and uses it as the decision variable. When the product's lead time is over, the order quantity is added to the inventory. This is done for an entire year.

The processes were parallelized using delayed and parallel functions, allowing for the simultaneous execution of numerous simulations and a high iteration rate. By distributing the simulations across multiple CPU cores, the total execution time was reduced significantly. To emphasize extreme demand values, the simulation process used a mixture of two normal distributions, one centered around the historical mean and another around a higher demand value. This approach helped in capturing the tail behavior of the demand distribution. This was crucial to understand the risk of stockouts and high-demand scenarios. **Figure 4** displays the working simulation model which assumes varied demand and lead time.

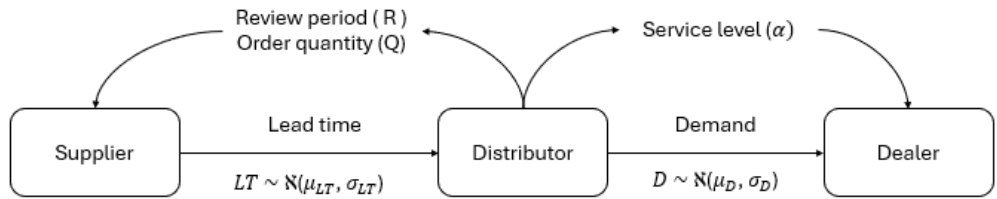

**Figure 4. Inventory model**

### 2.4.1 Constraints and Assumptions

Several assumptions are made in this study which include:

1) The demand for products follows a log-normal distribution.
2) The lead time for replenishment is deterministic and known in advance.
3) The model considers various costs associated with inventory management, including purchase costs, holding costs, and ordering costs. These costs are assumed to be known and constant over time, without considering factors such as inflation or economies of scale.
4) The model assumes that stockouts are not allowed, meaning that demand is always met without delay.
5) The demand for each product is assumed to be independent and identically distributed across periods.
6) The demand and the lead time are continuous random variables.

### 2.4.2. Convergence analysis

The below plots (**Figure 5**, **Figure 6**, **Figure 7**, and **Figure 8**) collectively help in assessing the convergence and reliability of MCS results. They provide insights into how quickly the estimates stabilize, the presence of biases or trends in the simulations, and the precision of the profit estimates.

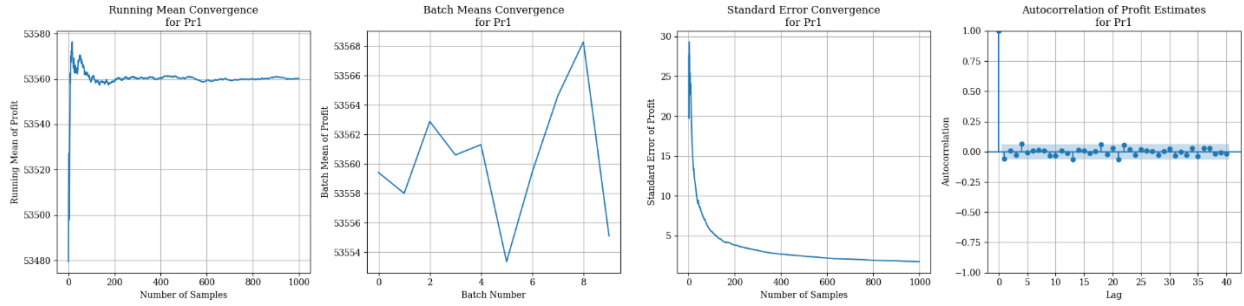

**Figure 5. Convergence plots Pr1**

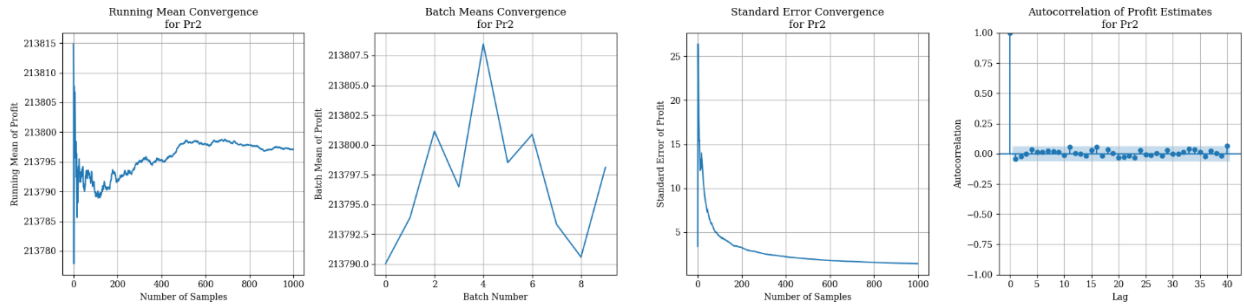

**Figure 6. Convergence plots Pr2**

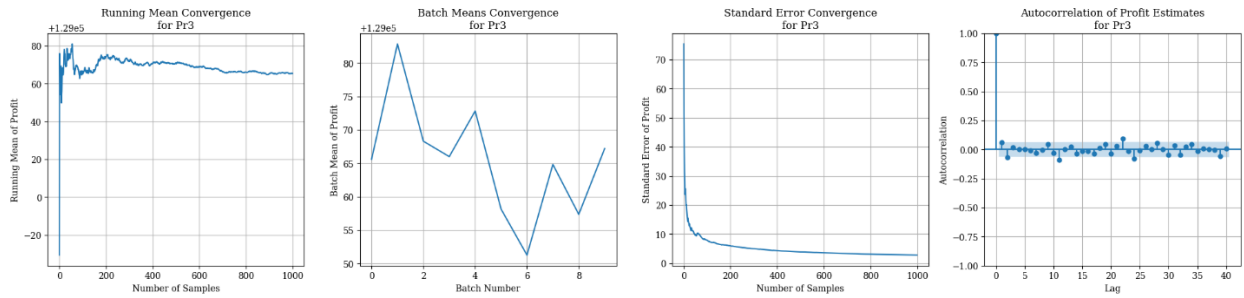

**Figure 7. Convergence plots Pr3**

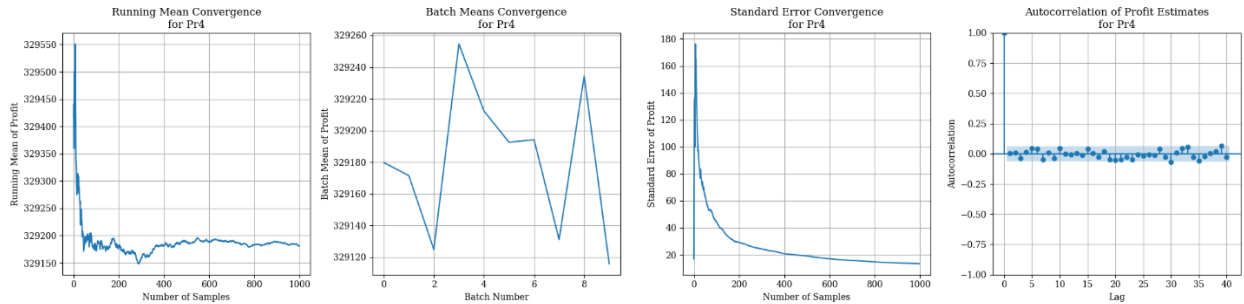

**Figure 8. Convergence plots Pr4**

The plots show that initially, the running mean fluctuates significantly, indicating high variability in early samples. As the number of samples increases, the running mean stabilizes, suggesting that the estimate of the average profit is converging. This indicates that the simulation is reaching a steady state where the mean profit estimate becomes reliable.

The batch means convergence plots show the mean profit for different batches of simulations. Large fluctuations might indicate variability in batches or the need for a larger batch size to smooth out the results.

The standard error plots show the standard error decreases as the number of samples increases, indicating that the estimate of the profit is becoming more precise with more simulations.

The autocorrelation measures how the profit estimates are correlated with themselves over different time lags. Most of the points are close to zero, indicating that the profit estimates are largely independent over time. This suggests that there is no significant time-based pattern or trend in the profit estimates, which is desirable in MCS as it indicates that each sample is independent.

## 3. Empirical analysis

Besides random sampling, this study experimented with conditional sampling to highlight critical periods of high demand. It starts with the expected demand and updates it based on the actual demand observed on the day the order is triggered. If an order is triggered, the demand for subsequent days is set to last high demand. If an order is not triggered, the model samples daily demand using a normal distribution. **Table 4** reports the output from both sampling approaches using the (p, Q) policy.

Table 4. Periodic review - Random vs. Conditional sampling for stochastic demand.

| Product | Reorder quantity | $\mu_{profit}$ | $\sigma_{profit}$ | $D_{LT}$ |
|---|---|---|---|---|
| | **Random Sampling** | | | |
| Pr1 | 4120 | 189543.97 | 4927.64 | |
| Pr2 | 33730 | 577530.97 | 5935.75 | |
| Pr3 | 7770 | 422764.87 | 22731.94 | |
| Pr4 | 2010 | 514094.89 | 72598.22 | |
| | Total profit → 1,703,934.7 | Time to run: 4527.05 secs | | |
| | **Conditional Sampling** | | | |
| Pr1 | 4120 | 189694.80 | 4802.71 | |
| Pr2 | 33740 | 577356.19 | 2646.94 | |
| Pr3 | 7790 | 420784.23 | 24457.91 | |
| Pr4 | 2080 | 429674.36 | 86558.09 | |
| | Total profit → 1,615,108.92 | Time to run: 3656.68 secs | | |

While in execution, the conditional sampling is much faster (3656.68 secs) compared to random sampling (4527.05 secs); however, the total expected profit is lower in the case of conditional sampling compared to random sampling (1,615,108.92 vs 1,703,934.7). This suggests that under the conditions simulated, conditional sampling may lead to slightly lower profits overall.

In the (r, Q) logic the reorder point (ROP) is estimated as $ROP = SS + D_{LT}$, where $D_{LT}$ = demand during lead time. The inventory level is checked daily, if the $INV \leq ROP$ and no order in hand, an order is triggered. Q (order quantity) is calculated to bring the inventory back up to the SS level plus the expected $D_{LT}$ to ensure a minimum order quantity. So, the (r, Q) system ensures that inventory levels are managed dynamically based on daily demand and lead time considerations, allowing for effective inventory control while meeting customer demand and minimizing costs. **Table 5** reports the output from (r, Q) policy.

**Table 5. Continuous review for stochastic demand.**

| Product | r | $\mu_{profit}$ | $\sigma_{profit}$ |
|---|---|---|---|
| Pr1 | 4095 | 447626.66 | 20366.69 |
| Pr2 | 33940 | 1633884.81 | 53635.00 |
| Pr3 | 7710 | 504179.23 | 4608.83 |
| Pr4 | 2370 | 688613.93 | 122741.87 |
| Total profit → 3,274,304.64 | Time to run: 4693.92 secs ||||

The output indicates a significant increase in expected profit (3,274,304.64) compared to previous simulations. This improvement suggests that for the given business case, the (r, Q) policy is effective in optimizing inventory decisions.

**Figure 9** supports the output of **Table 5** showing the relationship between order quantity (Q), reorder point (r), and expected profit for Pr1. It shows a smooth transition of expected profit across different combinations, with the highest profit in lighter areas and lower profits in darker colors. The optimal combination is in the middle of the plot, suggesting ordering quantities around 3000-3250 and reorder points around 4200-4250 are optimal for maximizing profit. The plot indicates sensitivity to changes in Q and r, suggesting a degree of robustness in the inventory policy. This information can help identify the best order quantity and reorder point to maximize profit, reducing the risk of significant profit loss due to minor adjustments.

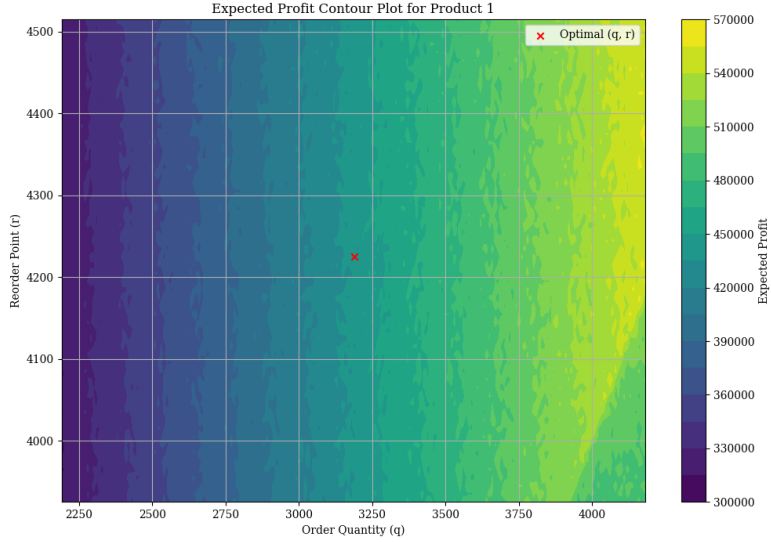

**Figure 9. Contour plot of expected profit (Pr1).**

A comparative analysis of optimization approaches by employing grid-search and Bayesian optimization was performed. Grid-search was employed over a range of Q and r values (incremented by 10). The objective was to evaluate the profitability and operational metrics (e.g., $\mu_{profit}$, $\sigma_{profit}$, $\mu_{SS}$, and mean proportion of lost orders) for different combinations of Q (order quantity) and r (reorder point). Empirical studies have found Bayesian optimization is a powerful approach for the global derivative-free optimization of non-convex expensive functions [45]. **Table 6** displays the consolidated output.

**Table 6. Performance optimization.**

| Products | $\mu_{profit}$ | r | Q | $\mu_{SS}$ | $\sigma_{profit}$ |
|---|---|---|---|---|---|
| | | **Grid search** | | | |
| Pr1 | 455298.75 | 500 | 3200 | 422 | 18893.89 |
| Pr2 | 2028436.29 | 1000 | 32500 | 352 | 4024.65 |
| Pr3 | 521983.34 | 500 | 7160 | 359 | 13082.08 |
| Pr4 | 860626.65 | 500 | 1620 | 465 | 85071.92 |
| | | **Total Profit: 3,866,345.03** | | | |
| | | **Bayesian optimization** | | | |
| Pr1 | 521836.49 | 4439 | 3845 | 3687 | 34849.75 |
| Pr2 | 1958021.72 | 35000 | 25000 | 27906 | 109806.55 |
| Pr3 | 771964.35 | 7956 | 10000 | 6647 | 48815.27 |
| Pr4 | 632979.93 | 2848 | 1912 | 2782 | 125500.15 |
| | | **Total profit: 3,884,802.49** | | | |

To this end, this study implemented and evaluated two inventory management policies, (p, Q) and (r, Q), using Monte Carlo Simulation (MCS) combined with random and conditional sampling techniques. The (p, Q) policy, which involves periodic reviews and fixed intervals for placing orders, demonstrated that while conditional sampling is faster, random sampling yields slightly higher profits. In contrast, the (r, Q) policy, which dynamically places orders based on inventory levels, significantly outperformed the (p, Q) policy in terms of expected profit. This suggests that continuous review systems (r, Q) are more effective for managing inventory under stochastic demand. Additionally, a comparative analysis of optimization approaches revealed that Bayesian optimization marginally outperforms grid search, further enhancing profit margins. However, given the stochastic nature of the algorithm, grid search may bring a similar level of improvement. Grid search is easily applicable and clearly explainable and provides a straightforward approach for identifying optimal reorder quantities and points. It is particularly useful in settings where interpretability and simplicity are prioritized.

Since MCS is quite versatile and popular at the same time in real-life business environments ([46]; [16]), a wide range of possible outcomes can be predicted by changing the assumptions and constraints of the proposed mode across all parameters. The model can be tailored to different scenarios and compared to numerous potential outcomes.

## 4. Limitations and Future Direction

There are still certain unanswered questions about simulation-based optimization techniques for stochastic inventory control. This proposed approach assumes that daily operations and decisions are based solely on the available inventory at that moment, without considering future demand forecasts or other complicated factors. Since numerous replications are needed to eliminate noise in the given result, the simulation is in general computationally expensive. Furthermore, simulation offers a "what-if" response to system inputs in which gradient information is not readily available, in contrast to mathematical models. Moreover, the simulation assumed a log-normal demand distribution, relying on historical data and static lead times for inventory replenishment. The analysis also overlooked costs such as storage, handling, and stockout penalties. External factors like market dynamics, supplier variability, simulation accuracy, historical data quality, and policy variability are also overlooked in this study. While the study provides valuable insights but requires careful interpretation and potential adjustments when applying the findings to real-world scenarios, The simulation typically works with one or more objective functions, several constraints, and numerous parameters. However, it becomes impractical to analyze every potential solution when there are a lot of parameter combinations or when one or more parameters have an excessive number of values [47].

**Conclusion**

The study explores the optimization of inventory management strategies using Monte Carlo Simulation (MCS) and applying both (p, Q) and (r, Q) inventory policies. The goal was to determine the optimal reorder quantities (Q) that maximize profit while minimizing stockout risks and excessive inventory holdings. The results showed that the (r, Q) policy significantly increased expected profit compared to the (p, Q) policy. This is due to its dynamic adaptation to fluctuating demand, resulting in better inventory control and higher profitability. The study also found that conditional sampling in MCS showed faster execution times but slightly lower overall profits. However, the study acknowledges limitations, including assumptions of deterministic lead times, log-normal demand distributions, and the exclusion of certain costs. Future research should incorporate more complex factors and explore other advanced optimization techniques. It demonstrates that dynamic inventory management can reduce stockout risks and excess inventory, leading to higher profitability. Monte Carlo Simulation (MCS) allows companies to model and analyze inventory strategies under different demand scenarios, providing a robust framework for informed decisions. Optimizing reorder quantities and reorder points can reduce costs associated with excess inventory and mitigate losses from stockouts. The study also emphasizes risk management, ensuring businesses are better equipped to handle high-demand scenarios without compromising service levels. The open-source nature of the framework makes it suitable for various industries and business sizes.

Data and source code: https://github.com/saritmaitra/MCS-Inventory-Management